\documentstyle[12pt]{article}

\author{Germano D'Abramo
\thanks{I am very grateful to my longlasting friend Eng.~Giorgio Borghini.}\\
\small IASF - Istituto Nazionale di Astrofisica,\\
\small Via Fosso del Cavaliere 100,\\
\small 00133 Roma, Italy}
\title{\bf The Seven Messengers and the ``Buzzati sequence''}
\date{}

\begin{document}

\maketitle
\begin{abstract}

A young Prince decides to explore the Father's Kingdom and aims to 
reach its furthermost boundaries. He starts from the City with a Caravan 
and seven fast and strong Messengers. They have the task to maintain the 
communications between the Caravan and the City during the exploration, 
going back and forth between the Caravan and the City while the Caravan 
inexorably goes away. 

Drawing inspiration from a fantasy tale by the Italian writer Dino 
Buzzati, I derive the geometric progression (named ``Buzzati sequence''
in his honor) which governs the duration of the Messenger's trips, 
renewing further the fascination of the tale. I also note with wonder how 
all this apparently hidden mathematical structure was already known to 
the author. 
An extension of the ``Buzzati sequence'' to relativistic velocities of the 
Caravan and the Messengers is finally presented as exercise.

\end{abstract}


\section{Introduction}

Some years ago a dear friend of mine gave me a collection of fantasy tales 
by Dino Buzzati, named `{\it La Boutique del Mistero}' \cite{Buzz}. Dino 
Buzzati (San Pellegrino {\it 1906-- }Milano{\it \ 1972}) was an italian 
journalist, a writer, a poet and a painter. He had been an extraordinary 
explorer of the human mind, of his anguish and fears, and his writings 
perfectly reflect his deep vision of the human condition. I already had 
the pleasure to read something by Buzzati and to go into his troubled 
and fascinating {\it Weltanschauung}, but the first tale of this 
collection, titled `{\it I sette messaggeri' } (The seven messengers) 
strongly impressed my `mathematical imagination': just more than thirty 
years old, a Prince decides to explore the Father's Kingdom and aims to 
reach its furthermost boundaries. He starts from the City with a Caravan 
and seven fast and strong Knights. These Knights have the task to 
maintain the communications between the Caravan and the City during the 
exploration, in other words they are the Messengers. The first Messenger 
starts two days after the departure of the Caravan, the second three 
days after and so on. Many years pass and the Kingdom seems to be 
endless. Even though the Messengers ride night and day, every day, their 
encounters with the Caravan become ever more rare and the contacts 
between the Caravan and the City absurdly far away in time.

Soon after having read the tale I tried to derive the time elapsed, for 
each Messenger, between each tour using the informations given by the 
author and, to my amazement, I obtained really the same values described 
in the tale. My wonder was mainly due to the fact that, by my 
experience, this particular attention to mathematical accuracy is a bit 
unusual in this kind of literture.

In the following Section I derive the mathematical expression of the 
duration of the tours. Let me refer to this relation as to the ``Buzzati 
sequence'', named after the great writer. Suggested by the fascinating 
atmosphere of the tale, where men seem to be lost in time, in the third 
Section I make the exercise of extending the previous formula to cases in 
which the velocities of the Caravan and Messengers are relativistic (namely,
closer to the velocity of light).

\section{Buzzati sequence}

Before deriving the ``Buzzati sequence'' it is appropriate to clarify all 
the approximations introduced with the aim of simplifying the 
derivation. I assume that all the velocities are constant in modulus 
($V_c$ is the velocity of the Caravan, while $V_m$ is the velocity of 
all Messengers. Obviously it must be $V_m>V_c$, otherwise the Caravan 
misses the Messengers) and that Messengers change direction {\bf 
instantaneously} soon after having reached the City or the Caravan (in 
other words they never stop). With $T_{n,i}$ I represent the time when 
the {\it i--}th Messenger leaves the Caravan for the {\it n}--th tour 
(for sake of simplicity time counting starts when the Caravan{\sl \ 
}leaves the City) and $\Delta T_{n,i}$ is the duration of the {\it 
n}--th tour (i.e. $\Delta T_{n,i}=T_{n+1,i}-T_{n,i}$).

With this assumptions it is not difficult to see that for first tours 
the following relation holds

\begin{equation}
\label{uno}_{}V_m\Delta T_{1,i}=2V_cT_{1,i}+V_c\Delta T_{1,i}, 
\end{equation}
which simply states that the distance covered by the {\it i--}th 
Messenger in his first tour must be equal to two times the distance 
covered by the Caravan just before this departure plus the distance 
covered by the Caravan during Messenger's travel. So from equation 
(\ref{uno}) we have for $\Delta T_{1,i}$

\begin{equation}
\label{due}\Delta T_{1,i}=\frac{2V_cT_{1,i}}{V_m-V_c}=
2qT_{1,i},{\rm if\,\,we\,\,define\,\,}q\equiv \frac{V_c}{V_m-V_c}. 
\end{equation}

The general expression for $\Delta T_{n,i}$ can be easily obtained from 
equation (\ref{due}) by substituting $T_{1,i}$ with the time when the 
{\it n--}th departure takes place. Proceeding by steps

$$
\begin{array}{c}
\Delta T_{2,i}=2q(T_{1,i}+\Delta T_{1,i})=\Delta T_{1,i}+2q\Delta
T_{1,i}=\Delta T_{1,i}(1+2q), 
\end{array}
$$
since the {\it i--}th Messenger starts his second tour at
$T_{2,i}=T_{1,i}+\Delta T_{1,i}$,

\begin{eqnarray*}
\Delta T_{3,i} & = &2q(T_{1,i}+\Delta T_{1,i}+\Delta T_{2,i})=\Delta
T_{2,i}+2q\Delta T_{2,i}=\Delta T_{2,i}(1+2q)=\Delta T_{1,i}(1+2q)^2, \\ 
\Delta T_{4,i} &=& 2q(T_{1,i}+\Delta T_{1,i}+\Delta T_{2,i}+\Delta
T_{3,1})=\Delta T_{3,i}+2q\Delta T_{3,i}=\Delta T_{3,i}(1+2q)=\\ 
 &=&\Delta T_{1,i}(1+2q)^3, \\ 
...., & &
\end{eqnarray*}

and so on. So we can write by induction

\begin{equation}
\label{tre}\Delta T_{n,i}=\Delta T_{1,i}(1+2q)^{n-1}=2qT_{1,i}(1+2q)^{n-1}. 
\end{equation}

Moreover comparing equation (\ref{due}) with equation (\ref{tre}) it is 
easy to see that the time $T_{n,i}$ when the {\it n--}th departure takes 
place is simply
\begin{equation}
\label{quattro}T_{n,i}=T_{1,i}(1+2q)^{n-1}, 
\end{equation}
which is a straightforward geometric progression.

If we adopt the relation $V_m=3/2V_c$ (so $q=2$) and we put $%
T_{1,i}=(i+1)days,i=1,...,7,$ like in the tale by Buzzati, we can 
verify how fast (\ref{quattro}) grows with {\it n. }In the following 
table some values of $T_{n,i}$ are shown

\begin{center}
{\small 
\begin{tabular}{|l|l|l|l|l|l|l|l|} \hline
           $q=2$  &  {\sc Mess. 1}  &  {\sc Mess. 2} & {\sc Mess. 3} & {\sc Mess. 4} & {\sc Mess. 5} &{\sc Mess. 6} &{\sc Mess. 7} \\ \hline \hline
   { $T_{1,i}$}  & 2 days & 3 days& 4 days& 5 days & 6 days  & 7 days & 8 days\\ 
   { $T_{2,i}$} &  10 "    & 15 "    &  20 "  &  25 "   & 30 " & 35 " & 40 " \\
   { $T_{3,i}$} &  50 "    & 75 "    & 100 "  &  125 "  & 150 " & 175 "  & 200 "\\
   { $T_{4,i}$} &  250 "  & 375 "   &  $\sim$ 1.4 yrs  & $\sim$ 1.7 yrs  & $\sim$ 2.05 yrs &  $\sim$ 2.4 yrs & $\sim$ 2.7 yrs \\ 
   { $T_{5,i}$} &  $\sim$ 3.4 yrs &$\sim$ 5.1 yrs     & $\sim$ 6.8 "  &  $\sim$ 8.6 "   & $\sim$ 10.3 " & $\sim$ 12.0 " & $\sim$ 13.7 " \\ 
   { $T_{6,i}$} &  $\sim$ 17.1 "  &$\sim$  25.7 "     & $\sim$ 34.2 " & $\sim$ 42.8 "   & $\sim$ 51.4 " & $\sim$ 60.0 "& $\sim$ 68.5 "\\ 
   { $T_{7,i}$} &  $\sim$  85.6 " &$\sim$ 128.4 "     & $\sim$ 171.2 "  & $\sim$ 214.0 " & $\sim$ 256.8 "  & $\sim$ 299.6 " & $\sim$ 342.5 "\\ \hline \hline
\end{tabular}
}
\end{center}

\section{Relativistic treatment}

It is appropriate to note that the way in which relations (\ref{tre}) 
and (\ref{quattro}) have been derived in the previous section doesn't 
take differences between `reference systems' into account, i.e. the time 
$T_{n,i}$ is the same for the City, the Caravan{\sl \ }and the {\it 
i--}th Messenger; moreover (\ref{tre}) and (\ref{quattro}) seem to hold 
for every absolute value of $V_m$ and $V_c$.

However, form the Theory of Special Relativity (TSR, A. Einstein, 1905) it 
is well known that nothing can overcome the speed of light and that the 
measure of space and time intervals depends on the particular reference 
system.

Here we are dealing with special relativity so time and space 
intervals will be measured with respect to inertial reference systems, 
that is, systems in which free motion of bodies (motion not subject to 
any kind of external force) is a constant speed motion.

According to the theory what is really invariant is the infinitesimal
`length' element in space--time 
\begin{equation}
\label{cinque}ds^2=c^2dt^2-dx^2-dy^2-dz^2, 
\end{equation}
namely, the `distance' between two events infinitesimally close in space 
and time. Indeed, in any inertial reference system we can assign a time 
and three spatial coordinates ($t,x,y,z$) to any event and (\ref{cinque}) 
can be written for two events infinitesimally close; the main 
point is that the value of $ds^2$ is the same for all those systems. So, 
for two inertial reference systems, say $O$ and $\overline{O}$, the following 
holds
\begin{equation}
\label{sei}c^2dt^2-dx^2-dy^2-dz^2=c^2\overline{dt}^2-\overline{dx}^2- 
\overline{dy}^2-\overline{dz}^2, 
\end{equation}
no matter what is the constant relative velocity between $O$ and 
$\overline{O}$.

From (\ref{sei}) it is not difficult to obtain the expression for the 
`proper time' of an observer, i.e. the time measured by a clock at rest 
in the observer reference system (for example, like we will see, the 
proper time of the Caravan or Messengers). In our inertial system ($O$) 
the infinitesimal displacement of the observer in space--time is of the 
form of (\ref{cinque}), while in the reference system of the observer 
($\overline{O}$), by definition, there is no spatial displacement and 
(\ref{cinque}) becomes
$$
\overline{ds}^2=c^2\overline{dt}^2. 
$$
But now we know that $ds^2=\overline{ds}^2$, so 
\begin{equation}
\label{sette}c^2\overline{dt}^2=c^2dt^2-dx^2-dy^2-dz^2. 
\end{equation}

Moreover if $\stackrel{\rightarrow}{v}(t)$ is the observer velocity
(not necessarily uniform) that we measure in our inertial system, we can 
write, for infinitesimal spatial displacements,
$$
dx=v_x(t)dt,\quad dy=v_y(t)dt,\quad dz=v_z(t)dt, 
$$
so, for proper time $\overline{t}$, we have (from (\ref{sette})) 
\begin{equation}
\label{otto}\overline{t}=\int\limits_0^t{\rm d}\overline{t}=\int\limits_0^t 
\sqrt{1-\frac{\stackrel{\rightarrow }{v}(\tau )\cdot \stackrel{%
\rightarrow }{v}(\tau )}{c^2}}{\rm d}\tau . 
\end{equation}

For a more detailed description of these basic principles and a more 
rigorous derivation of the formulas you can see, for example, the book 
by L. D. Landau and E. M. Lif\v sits, {\it The Classical Theory of 
Fields} \cite{Landau}.

Now we are able to deal with the extension of `Buzzati sequence' to 
cases in which $V_c/c\rightarrow 1,$ $V_m/c\rightarrow 1$, where $c$ is 
the speed of light; I will use the same notations and approximations of 
Section 2.

Suppose you live in the City (our inertial reference system) and you 
want to find out the proper time of the Caravan and Messengers, then all 
you need is (\ref{tre}), (\ref{quattro}) and (\ref{otto}). It is natural 
to consider $T_{1,i}$ like a time measured in the reference frame of 
the Caravan (in fact before their first tours the Messengers belong to 
the Caravan), while in the City we have, by (\ref{otto})
\begin{equation}
T_{1,i}=\sqrt{1-\frac{V_{c^{}}^2}{c^2}}T_{1,i,{\rm City}}\quad
\rightarrow \quad T_{1,i,{\rm City}}=
\frac 1{\sqrt{1-\frac{V_{c^{}}^2}{c^2}}}T_{1,i}. 
\label{nove}
\end{equation}

Obviously I suppose that all the clocks have been synchronized at the origin.
So for $T_{n,i,{\rm City}}$, which gives the time of the {\it n--}th
departure in the reference system of the City, we have 
\begin{equation}
\label{dieci}T_{n,i,{\rm City}}=T_{1,i,{\rm City}}(1+2q)^{n-1}=\frac 1{\sqrt{
1-\frac{V_{c^{}}^2}{c^2}}}T_{1,i}(1+2q)^{n-1}. 
\end{equation}

Multiplying (\ref{dieci}) by $\sqrt{1-V_c^2/c^2}$ we will obtain the proper
time $T_{n,i}$ of the Caravan which is obviously $T_{1,i}(1+2q)^{n-1}$ again$%
.$ Concerning the proper time of Messengers the situation is quite
different. The proper time $T_{n,i,{\rm Mes}}$ when Messengers start
their {\it n--}th tour, is the sum of the following two terms (see (\ref
{otto}) again) 
\begin{equation}
\label{undici}T_{n,i,{\rm Mes}}=\int\limits_0^{T_{1,i,{\rm City}}}\sqrt{1- 
\frac{V_{c^{}}^2}{c^2}}{\rm d}\tau \quad +\quad \int\limits_{T_{1,i,{\rm City%
}}}^{T_{n,i,{\rm City}}}\sqrt{1-\frac{V_m^2}{c^2}}{\rm d}\tau , 
\end{equation}
since before $T_{1,i,{\rm City}}$ the Caravan and the Messengers travel
together at the same speed $V_c$ . The integration (\ref{undici}) is trivial
and using equations (\ref{nove}) and (\ref{dieci}) we obtain 
\begin{equation}
\label{dodici}T_{n,i,{\rm Mes}}=\frac{\sqrt{1-\frac{V_m^2}{c^2}}}{\sqrt{1- 
\frac{V_{c^{}}^2}{c^2}}}T_{1,i}(1+2q)^{n-1}\quad +\quad \left( 1-\frac{\sqrt{
1-\frac{V_m^2}{c^2}}}{\sqrt{1-\frac{V_{c^{}}^2}{c^2}}}\right) T_{1,i}. 
\end{equation}

Now it is easy to verify that both (\ref{dieci}) and (\ref{dodici}) reduce
themselves to (\ref{quattro}) when $V_c/c\rightarrow 0$ and 
$V_m/c\rightarrow 0.$ Moreover the time $T_{n,i,{\rm Mes}}$ is asymptotically
lower than both proper time of the Caravan and $T_{n,i,{\rm City}}$ (this
phenomenon is of the same kind of the so--called `Twin Paradox').

A little more difficult exercise is to derive the expression of proper time
of Messengers when they reach the City during their {\it n--}th tour ($%
T_{n,i,{\rm Mes}}^{{\rm City}}$) . In the reference system of the City the 
{\it i--}th Messenger starts his {\it n--}th tour at $T_{n,i,{\rm City}}$,
when the Caravan is at distance $T_{n,i,{\rm City}}\times V_c$, so he will
reach the City after $T_{n,i,{\rm City}}\times V_c/V_m$ or $T_{n,i,{\rm City}%
}\times q/(1+q)$ (see (\ref{due})). Thus the proper time will be, using (\ref
{otto}) again,
$$
T_{n,i,{\rm Mes}}^{{\rm City}}=T_{n,i,{\rm Mes}}\quad +\quad \sqrt{1-\frac{
V_m^2}{c^2}}\frac q{1+q}T_{n,i,{\rm City}}, 
$$
which becomes, substituting $T_{n,i,{\rm City}}$ and $T_{n,i,{\rm Mes}}$with
(\ref{dieci}) and (\ref{dodici}), 
\begin{equation}
\label{tredici}T_{n,i,{\rm Mes}}^{{\rm City}}=\frac{\sqrt{1-\frac{V_m^2}{c^2}
}}{\sqrt{1-\frac{V_{c^{}}^2}{c^2}}}\left( 1+\frac q{1+q}\right)
T_{1,i}(1+2q)^{n-1}\quad +\quad \left( 1-\frac{\sqrt{1-\frac{V_m^2}{c^2}}}{
\sqrt{1-\frac{V_{c^{}}^2}{c^2}}}\right) T_{1,i}.\quad 
\end{equation}

In the following table some values of $T_{n,4,{\rm City}}$, $T_{n,4,{\rm Mes}%
}$, $T_{n,4}$ and $T_{n,4,{\rm Mes}}^{{\rm City}}$ for Messenger $4$ are shown; 
$q$ is still equal to $2$ and for $V_c$ I have chosen a test value of $c/2$ 
(so $V_m=3/4c$)

\begin{center}
{\small 
\begin{tabular}{|l|l|l|l|l|} \hline
$n$  &  $T_{n,4,{\rm City}}$  &  $T_{n,4,{\rm Mes}}$ &  $T_{n,4}$ (non relativistic)&  $T_{n,4,{\rm Mes}}^{{\rm City}}$ \\ \hline \hline
   1  & 5.8 days & 5 days& 5 days& 7.54 days \\ 
   2  &  28.9 "    & 20.3 "    & 25  "  &  33.0 "    \\
   3  &  144.3 "    & 96.6 "    & 125 "  & 160.3 "  \\
   4  &  $\sim$ 1.9 yrs  & $\sim$ 1.3 yrs  &  $\sim$ 1.7 yrs  & $\sim$ 2.2 yrs    \\ 
   5  &  $\sim$ 9.9 " &$\sim$ 6.5 "     & $\sim$ 8.6 "  &  $\sim$ 10.9 "    \\ 
   6  &  $\sim$ 49.4 "  &$\sim$ 32.7 "     & $\sim$ 42.8 " & $\sim$ 54.5 "  \\ 
   7  &  $\sim$ 247.1 " &$\sim$ 163.5 "     & $\sim$ 214.0 "  & $\sim$ 272.5 " \\ \hline \hline
\end{tabular}
}
\end{center}

Lastly suppose that the Caravan uses electromagnetic waves to 
communicate, namely $V_m=c$. Hence, it is easy to see that, at first 
order in $V_c/c$, equation (\ref{dieci}) becomes $T_{n,i,{\rm 
City}}=T_{n,i}=T_{1,i}(1+2(n-1)V_c/c)$ (remember that the exchange of 
informations between Caravan and Messengers --electromagnetic Messengers 
too-- was supposed instantaneous). 

And the last thought, why not to see Buzzati's marvelous tale as a 
metaphor for the future space exploration of the mankind?

\end{document}